\documentclass[11pt]{article}
\usepackage{latexsym}
\usepackage{amsmath}
\usepackage{amssymb}
\usepackage{amscd}
\usepackage{bm}
\usepackage[dvips]{graphicx}

\newtheorem{Th}{Theorem}[section]
\newtheorem{Cor}{Corollary}[section]
\newtheorem{Prop}{Proposition}[section]
\newtheorem{Lem}{Lemma}[section]
\newcounter{Remark}[section]
\newcounter{Example}[section]
\newcommand{\qed}{\hfill $\Box$ \par}
\newcommand{\bet}{\begin{Th}}
\newcommand{\ent}{\stepcounter{Cor}
   \stepcounter{Prop}\stepcounter{Lem}
   \stepcounter{Remark}\stepcounter{Example}\end{Th}}

\newcommand{\bec}{\begin{Cor}}
\newcommand{\enc}{\stepcounter{Th}
   \stepcounter{Prop}\stepcounter{Lem}
   \stepcounter{Remark}\stepcounter{Example}\end{Cor}}
\newcommand{\bep}{\begin{Prop}}
\newcommand{\enp}{\stepcounter{Th}
   \stepcounter{Cor}\stepcounter{Lem}
   \stepcounter{Remark}\stepcounter{Example}\end{Prop}}
\newcommand{\bel}{\begin{Lem}}
\newcommand{\enl}{\stepcounter{Th}
   \stepcounter{Cor}\stepcounter{Prop}
   \stepcounter{Remark}\stepcounter{Example}\end{Lem}}

\newcommand{\Remark}{
   \stepcounter{Remark}
   \noindent{Remark\,  \thesection.\theRemark \, }:
   \stepcounter{Th}\stepcounter{Cor}\stepcounter{Prop}
   \stepcounter{Lem}\stepcounter{Example}}

\newcommand{\Proof}{\noindent{Proof\,}:\ }

\begin{document}
\def\R{{\mathbb{R}}}
\def\Si{\Sigma}
\def\si{\sigma}
\def\de{\delta}
\def\vp{\varphi}
\def\om{\omega}
\def\al{\alpha}
\def\be{\beta}
\def\si{\sigma}
\def\Si{\Sigma}
\def\pa{\partial}
\def\tr{\mbox{\rm tr}}
\def\bi{\bar i}
\def\ga{\gamma}
\def\Ga{\Gamma}
\def\eps{\epsilon}
\def\Om{\Omega}
\def\la{\lambda}
\def\La{\Lambda}
\def\de{\delta}
\def\De{\Delta}
\def\vp{\varphi}
\def\dps{\displaystyle}
\def\triad{\triangledown}
\def\ad{\mbox{\rm ad}}
\def\Ad{\mbox{\rm Ad}}
\def\tr{\mbox{\rm tr}}
\def\rank{\mbox{\rm rank}}
\def\Ker{\mbox{\rm Ker}}
\def\Tr{\mbox{\rm Tr}}
\def\RC{R_{\al \be \ga \de}}
\def\Ric{\mbox{\rm Ric}}
\def\bRic{\overline{\Ric}}
\def\ve{\varepsilon}
\def\e{\epsilon}
\def\s3{\sqrt3}
\def\dm{\mbox{\rm diag}}
\def\ve{\varepsilon}
\def\bpm{\begin{pmatrix}}
\def\epm{\end{pmatrix}}
\def\beq{\begin{equation}}
\def\enq{\end{equation}}
\def\bary{\begin{array}{ll}}
\def\eary{\end{array}}
\def\s{\sqrt}

%%%%%%%%%%%%%%

\title{The Dorfmeister-Neher theorem on isoparametric hypersurfaces}

\author{Reiko Miyaoka}

\maketitle
\abstract
{A new proof of the homogeneity of isoparametric 
hypersurfaces with six simple principal curvatures 
\cite{DN} is given in a method applicable to the 
multiplicity two case.}
\footnote{ 


Partially supported by Grants-in-Aid for Scientific Reseach, 
19204006, The Ministry of Education, Japan}

\section{Introduction}
\label{intro}
The classification problem of isoparametric hypersurfaces 
is remaining in some cases of four and six principal 
curvatures (see \cite{CCJ}, \cite{I}).
The homogeneity in the case $(g,m)=(6,1)$  
was proved by Dorfmeister-Neher \cite{DN}.  
A shorter proof was given in \cite{M2},  
but some argument was insufficient 
(pointed out by Xia Qiaoling).  
Moreover, we found it difficult to extend the   
method to the case $(g,m)=(6,2)$. 

In the present paper, we show that a 
 delicate change of signs of some vectors at 
 anti-podal points on a leaf, which is 
related to the back ground symmetry caused by a spin 
action, is essential. 
This investigation is also indispensable to attack 
on the case $m=2$. 
Before treating this overwhelmingly difficult case, 
a complete short proof for $m=1$ will give us an overview 
 how to settle the problem in the case $m=2$ \cite{M3}.

\S\ref{pre} $\sim$ \S\ref{kernel} consist of reviews of 
\cite{M1} and \cite{M2}. 
We do not repeat the proofs in \cite{M1}, but 
give those of \cite{M2} in a refined manner. 
The shape operators of each focal 
submanifold $M_{\pm}$  consist of an $S^1$-family of
isospectral transformations with simple eigenvalues 
$\pm\sqrt3$, $\pm 1/\sqrt3$, $0$. There are many such 
$S^1$-families  (see \S\ref{pre}), but in 
\S\ref{sdescription} $\sim$ \S\ref{sfinal},   
we narrow down them by using both local and global properties of isoparametric hypersurfaces, 
and conclude that non-homogeneous cases cannot occur.

%%%%%%%%%%%%
%%%%%%%%%%%%

\section{Preliminaries}
\label{pre}

We refer readers to \cite{Th} for a nice survey of 
isoparametric hypersurfaces. 
Here we review fundamental facts and the notation given in \cite{M1}.
Let $M$ be an isoparametric hypersurface in the 
unit sphere $S^{n+1}$, 
with a unit normal vector field $\xi$.
We denote the Riemmannian connection on 
$S^{n+1}$ by $\tilde\triangledown$, 
and that on $M$ by $\triangledown$. 
The principal curvatures of $M$ are given by constants 
$\lambda_1\ge \dots \ge \lambda_n$, 
and the curvature distribution for 
$\la\in \{\la_\al\}$ is denoted by $D_\la(p)$,  
 $m_\la=\dim D_\la(p)$.
In our situation,  
$D_\la$ is completely integrable and  
a leaf $L_\la$  of $D_\la$ is an 
$m_\la$-dimensional sphere of $S^{n+1}$. 
Choose a local orthonormal frame
$e_1,\dots,e_n$ consisting of unit principal vectors 
corresponding to $\la_1,\dots,\la_n$.
We express
\begin{equation}
\tilde \triangledown _{e_\al}e_{\beta}
=\Lambda_{\al \beta}^\sigma e_\si 
+\la_\al\de_{\al \beta}\xi,
\quad  \La_{\al \beta}^\ga = -\La_{\al \ga}^\beta,
\label{2.1}
\end{equation}
where $1\le\al,\beta,\si\le n$, using the Einstein convention. 
The curvature tensor $R_{\al \beta \ga \de}$ of $M$ is given by
\begin{equation}
\begin{array}{ll}
&\RC =(1+\la_\al\la_\beta)(\de_{\be \ga}\de_{\al \de}
-\de_{\al \ga}\de_{\be \de})\\
 & =e_\al(\La_{\be \ga}^\de)-e_\be(\La_{\al \ga}^\de)
+\La_{\be \ga}^\si\La_{\al\si}^\de-\La_{\al \ga}^\si\La_{\be \si}^\de
-\La_{\al \be}^\si\La_{\si \ga}^\de+\La_{\be \al}^\si\La_{\si \ga}^\de.
\end{array}
\label{2.2}
\end{equation}
From the equation of Coddazi, 
we obtain
\begin{equation}
e_\beta (\la_\al)=\La_{\al \al}^\beta(\la_\al -\la_\beta),
\qquad \text{for}
\quad \al\ne\beta, 
\label{2.5}
\end{equation}
and if $\la_\al, \la_\be, \la_\ga$ are distinct, 
we have
\begin{equation}
\La_{\al \be}^\ga(\la_\be-\la_\ga)=\La_{\ga \al}^\be(\la_\al-\la_\be)
=\La_{\be\ga}^\al(\la_\ga-\la_\al). 
\label{2.6}
\end{equation}
Moreover, 
\begin{equation}
\La_{ab}^\ga=0,\quad\La_{aa}^\ga=\La_{bb}^\ga,
\qquad \text{if}\quad \la_a
=\la_b\ne\la_\ga\quad\text{and}\quad a\ne b,  
\label{2.7}
\end{equation}
hold, and since $\la_\al$ is constant on $M$, it follows from (\ref{2.5}), 
\begin{equation}
\La_{\al \al}^\ga =0\qquad \text{if}\quad \la_\ga\ne \la_\al.
\label{2.8}
\end{equation}

\smallskip
When the number $g$ of principal curvatures is six, 
the multiplicity $m$ of $\la_i$ is independent of $i$  
and takes values 1 or 2 \cite{A}. 
In the following, let $(g,m)=(6,1)$. As is well known,
$\la_i=\cot(\theta_1+\displaystyle{\frac{(i-1)\pi}{6}}),\, 
1\le  i\le 6,\, 0<\theta_1<\frac{\pi}{6}$, modulo $\pi$. 
Since the homogeneity is independent of the choice of
$\theta_1$, we take 
$$
\theta_1=\dps\frac {\pi}{12}=-\theta_6,\quad  
\theta_2=\dps\frac {\pi}{4}=-\theta_5,\quad 
\theta_3=\dps\frac {5\pi}{12}=-\theta_4
$$
 so that
\begin{equation}
\la_1=-\la_6=2+\sqrt3,\quad \la_2=-\la_5=1,\quad
\la_3=-\la_4=2-\sqrt3.
\label{pc}
\end{equation}
Note that we choose 
$\theta_i\in (-\frac{\pi}{2},\frac{\pi}{2})$. 
By (\ref{2.7}) and (\ref{2.8}), a leaf $L_i=L_i(p)$ of $D_i(p)=D_{\la_i}(p)$ is a 
geodesic of the corresponding curvature sphere.

For $a=6$ or $1$, define the  focal map 
$f_a\colon M \rightarrow S^{7}$ by
$$
f_a(p)=\cos\theta_a p+\sin\theta_a \xi_p,
$$
which collapses $L_a(p)$ into a point $\bar p=f_a(p)$. 
Then we have
\begin{equation}
df_a(e_j) =\sin\theta_a(\la_a-\la_j)e_j,
\label{2.9}
\end{equation}
where the right hand side is considered as a vector in 
$T_{\bar p}S^{7}$ by 
a parallel translation in $S^7$. We always use such 
identification. The rank of $f_a$ is constant and we 
obtain the focal submanifold $M_a$ of $M$:
$$
M_a=\lbrace \cos\theta_a p+\sin\theta_a \xi_p \mid p\in M\rbrace.
$$
By (\ref{2.9}), the tangent space of $M_a$ is given by 
$T_{\bar p}M_a=\oplus_{j\ne a}D_j(q)$ for 
any $q\in f_a^{-1}(\bar p)$. 
 An orthonormal basis of the normal space of 
$M_a$ at $\bar p$ is given by 
\begin{equation}
\eta_q=-\sin\theta_a q+\cos\theta_a \xi_q,\quad \zeta_q=e_a(q)  
\label{normal}
\end{equation}
for any $q\in L_a(p)=f_a^{-1}(\bar p)$. 

Now, the connection
$\bar\triangledown$ on $M_a$ is induced
from the  connection $\tilde \nabla$, that is
$$
{1\over{\sin\theta_a(\la_a-\la_j)}}
\tilde \nabla_{e_{j}}X
=\bar\triangledown_{e_{ j}}\tilde X
+\bar\triad_{e_{ j}}^{\bot}\tilde X,
\quad \la_j\ne \la_a, 
$$
where $X$ is a tangent field on $S^7$  
in a neighborhood of $p$, and $\tilde X$ is
the one 
near $\bar p$ translated from $X$.
Note that 
$\bar\triad_{e_{ j}}^{\bot}\tilde X$ denotes the 
normal component in $S^{7}$.
In particular, we have for $j\ne a$, 
\begin{align}
\bar\triad_{e_{ j}}\tilde e_{ k} &
=\frac1{\sin\theta_a(\la_a-\la_j)}
\sum_{\l\ne a}\La_{ jk}^{ l}e_{ l},
\label{shape1}\\
\bar\triad_{e_{ j}}^{\bot}\tilde e_{ k} &
=\frac1{\sin\theta_a(\la_a-\la_j)}
\lbrace\La_{ j k}^ae_a
+\sin\theta_a(1+\la_j\la_a)\de_{ j k}\eta_p\rbrace,
\label{shape}
\end{align}
using $\langle\la_j\xi_p-p,\eta_p\rangle
=\sin\theta_a(1+\la_j\la_a)$. 
In the following, we identify $\tilde e_k$ with $e_k$. 
Denote by $B_N$ the shape operator of $M_a$ 
with respect to the normal 
vector $N$. 
Then from (\ref{shape1}) and (\ref{shape}),  we obtain:
\bel \mbox{\rm\cite{M1}(Lemma 3.1)}  When we 
identify $T_{\bar p}M_a$ with 
$\oplus_{j=1}^5D_{a+j}(p)$ where the indices 
are modulo 6, the second fundamental tensors $B_{\eta_p}$ and 
$B_{\zeta_p}$ at $\bar p$ 
are given respectively by 
$$
B_{\eta_p}=
\begin{pmatrix}
\sqrt 3 & 0 & 0 & 0 & 0\\
 0  & \frac1{\sqrt 3} & 0 & 0 & 0\\
 0 &0 & 0 & 0 & 0\\
 0 & 0 &  0  & -\frac1{\sqrt 3} & 0\\
 0  & 0 & 0 &  0 & -\sqrt 3
\end{pmatrix},
$$
$$
B_{\zeta_p}=
\begin{pmatrix}
0 & b_{a+1\,a+2} & b_{a+1\,a+3} & b_{a+1\,a+4} & b_{a+1\,a+5}\\
b_{a+1\,a+2} & 0 & b_{a+2\,a+3} & b_{a+2\,a+4} & b_{a+2\,a+5}\\
b_{a+1\,a+3} & b_{a+2\,a+3} & 0 & b_{a+3\,a+4} & b_{a+3\,a+5}\\
b_{a+1\,a+4} & b_{a+2\,a+4} & b_{a+3\,a+4} & 0 & b_{a+4\,a+5}\\
b_{a+1\,a+5} & b_{a+2\,a+5} & b_{a+3\,a+5} & b_{a+4\,a+5} & 0
\end{pmatrix},
$$
where 
\beq
b_{jk}=\frac1{\sin\theta_a(\la_a-\la_j)}\La_{jk}^a
=\frac1{\sin\theta_a(\la_j-\la_a)}\La_{ja}^k, 
\quad a=6, 1
\label{bij}
\enq
\label{l2.1}
\enl

In fact, from (\ref{shape}) it follows $B_{\eta_p}(e_{ j})=\mu_j e_j$, where for $a$ is, say 6,  
\begin{equation}
\mu_j=\frac{1+\la_j\la_6}{\la_6-\la_j},\quad
\mu_1=\s3=-\mu_5,\,
\mu_2=1/\s3=-\mu_4,\,\mu_3=0, 
\label{beta}
\end{equation}
and $b_{jk}=b_{kj}$ follows from (\ref{2.6}). 
In the following, we denote $M_+=M_6$ and $M_-=M_1$. 
Note that  both are minimal. 
It is easy to see that any unit normal vector 
is written as $\eta_q$ in (\ref{normal}) for some 
$q\in L_6(p)$, and we have immediately:

\bel \mbox{\rm\cite{Mu}, \cite{M1}} The shape operators are isospectral, 
i.e., the eigenvalues of $B_N$ are 
$\pm \sqrt 3,\, \displaystyle{\pm\frac1{\sqrt 3}},\, 0$, 
for any unit normal $N$. 
\label{l2.2}
\enl

 For a fixed $p\in f_a^{-1}(\bar p)$,  
 all the shape operators for unit normals at $\bar p$ 
 are expressed as 
\begin{equation}
L(t)=\cos tB_{\eta_p}+\sin tB_{\zeta_p},\quad t\in [0,2\pi).
\label{Lt}
\end{equation}
The homogeneous hypersurfaces $M^h$ with $(g,m)=(6,1)$ 
are given as the principal orbits of the isotropy action 
 of the rank two symmetric space $G_2/SO(4)$, 
where two singular orbits correspond to the focal 
submanifolds $M^h_{\pm}$. 
In \cite{M1}, we show that the shape operators of $M^h_+$ and $M^h_-$ are given respectively by :
\beq
\begin{array}{ll}
\cos t
\begin{pmatrix}
\sqrt 3 & 0 & 0 & 0 & 0\\
 0  & \frac1{\sqrt 3} & 0 & 0 & 0\\
 0 &0 & 0 & 0 & 0\\
 0 & 0 &  0  & -\frac1{\sqrt 3} & 0\\
 0  & 0 & 0 &  0 & -\sqrt 3
\end{pmatrix}+\sin t\bpm
0&0&0&0&\sqrt3\\
0&0&0&\frac1{\sqrt3}&0\\
0&0&0&0&0\\
0&\frac1{\sqrt3}&0&0&0\\
\sqrt3&0&0&0&0
\epm,\\
\cos t
\begin{pmatrix}
\sqrt 3 & 0 & 0 & 0 & 0\\
 0  & \frac1{\sqrt 3} & 0 & 0 & 0\\
 0 &0 & 0 & 0 & 0\\
 0 & 0 &  0  & -\frac1{\sqrt 3} & 0\\
 0  & 0 & 0 &  0 & -\sqrt 3
\end{pmatrix}+\sin t\bpm
0&1&0&0&0\\
1&0&0&-\frac2{\sqrt3}&0\\
0&0&0&0&0\\
0&-\frac2{\sqrt3}&0&0&1\\
0&0&0&1&0
\epm. 
\end{array}
\label{dimE2B}
\enq
These imply that $M_{\pm}$ are {\em not} congruent to 
each other. 

Note that there {\it exist} many other one parameter 
families of isospectral operators $\cos t B_\eta+\sin tA$, 
where, for instance, $A$ is given by  
\beq
\begin{array}{ll}
\bpm
0&0&-\sqrt{\frac32}&0&0\\
0&0&0&\frac1{\sqrt3}&0\\
-\sqrt{\frac32}&0&0&0&\sqrt{\frac32}\\
0&\frac1{\sqrt3}&0&0&0\\
0&0&\sqrt{\frac32}&0&0
\epm,\quad
\bpm
0&0&0&0&\sqrt3\\
0&0&-\frac1{\sqrt6}&0&0\\
0&-\frac1{\sqrt6}&0&\frac1{\sqrt6}&0\\
0&0&\frac1{\sqrt6}&0&0\\
\sqrt3&0&0&0&0
\epm,\\
\bpm
0&\frac5{3\sqrt3}&0&\frac2{3\sqrt3}&0\\
\frac5{3\sqrt3}&0&\frac4{3\sqrt3}&0&-\frac2{3\sqrt3}\\
0&\frac4{3\sqrt3}&0&\frac4{3\sqrt3}&0\\
\frac2{3\sqrt3}&0&\frac4{3\sqrt3}&0&-\frac5{3\sqrt3}\\
0&-\frac2{3\sqrt3}&0&-\frac5{3\sqrt3}&0\epm,
\end{array}
\label{sample}
\enq
and so forth. We see in the homogeneous case, 
the kernel does not depend on $t$, while it depends in other cases. 
In the following, we show that all the latter cases are 
not admissible 
to the shape operators of the focal submanifolds of 
isoparametric hypersurfaces with $(g,m)=(6,1)$.

\section{Isospectral operators and Gauss equation}
\label{sisospectral}

 By Lemma \ref{l2.2},
$L(t)=\cos tB_\eta+\sin tB_\zeta$ is isospectral and so 
can be written as
\begin{equation}
L(t)=U(t)L(0)U^{-1}(t)\label{2.11}
\end{equation}
for some $U(t)\in O(5)$.
Moreover, this implies the Lax equation
\begin{equation}
L_t(t)=\frac{d}{dt}L(t)=[H(t),L(t)],\label{2.12}
\end{equation}
where
$$
H(t)=U_t(t)U(t)^{-1}\in{\frak o(5)}.
$$
In particular, we have $L(0)=B_\eta$, 
and
\begin{equation}
L_t(t)=-\sin tB_\eta+\cos tB_\zeta=L(t+\pi/2),
\label{L'}
\end{equation}
and hence for $L_t(0)=B_\zeta=(b_{ij})$, $b_{ij}=b_{ji}$,  and 
$H(0)=(h_{ij})$, $h_{ji}=-h_{ij}$, we can express 
\begin{equation}
\begin{array}{ll}
B_\zeta&=L_t(0)=[H(0),B_\eta]\\
&=\begin{pmatrix}
0&-\frac2{\sqrt3}h_{12}&-\sqrt3h_{13}&
-\frac4{\sqrt3}h_{14}&-2\sqrt3h_{15}\\
\frac2{\sqrt3}h_{21}&0&-\frac1{\sqrt3}h_{23}&
-\frac2{\sqrt3}h_{24}&-\frac4{\sqrt3}h_{25}\\
\sqrt3h_{31}&\frac1{\sqrt3}h_{32}&0&
-\frac1{\sqrt3}h_{34}&-\sqrt3h_{35}\\
\frac4{\sqrt3}h_{41}&\frac2{\sqrt3}h_{42}&
\frac1{\sqrt3}h_{43}&0&-\frac2{\sqrt3}h_{45}\\
2\sqrt3h_{51}&\frac4{\sqrt3}h_{52}&
\sqrt3h_{53}&\frac2{\sqrt3}h_{54}&0
\end{pmatrix}. 
\end{array}
\label{BH}
\end{equation}
Note that the eigenvectors of $L(t)$ are given by
\begin{equation}
e_{j}(t)=U(t)e_{j}(0),\label{2.14}
\end{equation}
which implies 
\beq
\nabla_{\frac{d}{dt}}e_{j}(t)=H(t)e_{j}(t).
\label{ddt}
\enq
Here we have 
\beq
\nabla_{\frac{d}{dt}}
=c_0 \nabla_{e_6},\quad c_0=\s2(\s3-1)/4, 
\label{c0}
\enq
because $L_6$ has radius $|\sin\theta_6|=c_0$. 
Hence we obtain 
\begin{equation}
H(0)=(c_0 \La_{6j}^{i}(0)),
\label{2.15}
\end{equation}
where $i$ denotes the row and $j$ denotes 
the column indices. 
Moreover, denoting the $(i,j)$ component of 
$L(t+\frac{\pi}{2})$ by $b_{ij}(t)$ where $b_{ji}(t)=b_{ij}(t)$, 
we have
$$
\bary
\left(\nabla_{e_6}L(t+\frac{\pi}{2})\right)_{ij}&=
e_6(b_{ij}(t))
-b_{kj}(t)\La_{6i}^k(t)-b_{ik}(t)\La_{6j}^k(t)\\
&=e_6(b_{ij}(t))
+\La_{6k}^{i}(t)b_{kj}(t)
-b_{ik}(t)\La_{6j}^{k}(t).
\eary
$$
Because  
$L_t(t+\frac{\pi}{2})=c_0\nabla_{e_6}L(t+\frac{\pi}{2})$,  
 $L_t(\pi/2)=-B_\eta$ and  
$L(\pi/2)=B_\zeta$, 
multiplying $-c_0$ to the both sides and putting 
$t=0$, we obtain
\begin{equation}
B_\eta
=-c_0e_6(B_\zeta)
-[H(0),B_\zeta].
\label{2.16}
\end{equation}
Now, rewrite (\ref{BH}) as 
$$
H(0)=\begin{pmatrix}
0&-\frac{\sqrt3}{2}b_{12}&-\frac1{\sqrt3}b_{13}&
-\frac{\sqrt3}{4}b_{14}&-\frac1{2\sqrt3}b_{15}\\
\frac{\sqrt3}{2}b_{21}&0 & -\sqrt3b_{23}& 
-\frac{\sqrt3}{2}b_{24}&-\frac{\sqrt3}{4}b_{25}\\
\frac1{\sqrt3}b_{31}&\sqrt3b_{32}&0& 
-\sqrt3b_{34}&-\frac1{\sqrt3}b_{35}\\
\frac{\sqrt3}{4}b_{41}&\frac{\sqrt3}{2}b_{42}
&\sqrt3b_{43}&0&-\frac{\sqrt3}{2}b_{45}\\
\frac1{2{\sqrt3}}b_{51}&\frac{\sqrt3}{4}b_{52}&
\frac1{\sqrt3}b_{53}&\frac{\sqrt3}{2}b_{54}&0
\end{pmatrix},
$$
 and substitute this into (\ref{2.16}). Then we have the following formulas which we use later :

[1.1]
$\sqrt3=2(\frac{\sqrt3}{2}b_{12}^2
+\frac1{\sqrt3}b_{13}^2
+\frac{\sqrt3}{4}b_{14}^2
+\frac1{2\sqrt3}b_{15}^2)$

[2.2] $\frac1{\sqrt3}=
2(-\frac{\sqrt3}{2}b_{21}^2
+\sqrt3b_{23}^2+\frac{\sqrt3}{2}b_{24}^2
+\frac{\sqrt3}{4}b_{25}^2)$

[3.3] $0=2(-\frac1{\sqrt3}b_{31}^2
-\sqrt3b_{32}^2+
\sqrt3b_{34}^2+
\frac1{\sqrt3}b_{35}^2)$

[4.4] $-\frac1{\sqrt3}=
2(-\frac{\sqrt3}{4}b_{41}^2
-\frac{\sqrt3}{2}b_{42}^2
-\sqrt3b_{43}^2
+\frac{\sqrt3}{2}b_{45}^2)$

[5.5]
$-\sqrt3=-2(\frac1{2\sqrt3}b_{51}^2
+\frac{\sqrt3}{4}b_{52}^2
+\frac1{\sqrt3}b_{53}^2
+\frac{\sqrt3}{2}b_{54}^2)$

[1.2]
$0=-c_0e_6(b_{12})
+\frac4{\sqrt3}b_{13}b_{32}
+\frac{3\sqrt3}{4}b_{14}b_{42}
+\frac5{4\sqrt3}b_{15}b_{52}$

[1.3]
$0=-c_0e_6(b_{13})
-\frac{\sqrt3}{2}b_{12}b_{23}
+\frac{5\sqrt3}{4}b_{14}b_{43}
+\frac{\sqrt3}{2}b_{15}b_{53}$

[1.4]
$0=-c_0e_6(b_{14})
-\frac2{\sqrt3}b_{13}b_{34}
+\frac2{\sqrt3}b_{15}b_{54}$

[1.5]
$ 0=-c_0e_6(b_{15})
+\frac{\sqrt3}{4}b_{12}b_{25}
-\frac{\sqrt3}{4}b_{14}b_{45}$

[2.3]
$ 0=-c_0e_6(b_{23})
-\frac5{2\sqrt3}b_{21}b_{13}
+\frac{3\sqrt3}{2}b_{24}b_{43}
+\frac7{\sqrt3}b_{25}b_{53}$

[2.4]
$0=-c_0e_6(b_{24})
-\frac{3\sqrt3}{4}b_{21}b_{14}
+\frac{3\sqrt3}{4}b_{25}b_{54}$

[2.5]
$0=-c_0e_6(b_{25})
-\frac2{\sqrt3}b_{21}b_{15}
+\frac2{\sqrt3}b_{23}b_{35}$

[3.4]
$0=-c_0e_6(b_{34})
-\frac{7}{4\sqrt3}b_{31}b_{14}
-\frac{3\sqrt3}{2}b_{32}b_{24}
+\frac{5}{2\sqrt3}b_{35}b_{54}$

[3.5]
$0=-c_0e_6(b_{35})
-\frac{\sqrt3}{2}b_{31}b_{15}
-\frac{5\sqrt3}{4}b_{32}b_{25}
+\frac{\sqrt3}{2}b_{34}b_{45}$

[4.5]
$0=-c_0e_6(b_{45})
-\frac{5}{4\sqrt3}b_{41}b_{15}
-\frac{3\sqrt3}{4}b_{42}b_{25}
-\frac4{\sqrt3}b_{43}b_{35}$

\noindent
These are nothing but another description of a part of 
the Gauss equations (\ref{2.2}). 
%(PENDING)
%The homogeneity follows if we have a local frame 
%with respect to which all $b_{ij}$ are constant, 
%which implies all $\La_{\al\be}^\ga$ are locally  
%constant. In that case, we can apply Singer's 
%strong curvature-homogeneous theorem \cite{KN}, p. , 
%and then the rigidity theorem on hypersurfaces. 
%We use this argument in \S\ref{sdimE2}. 
%%%%%%%%%%%%%%%%%%%%%%%%

\section{Global properties}
\label{GP}

An isoparametric hypersurface $M$ can be uniquely
extended to a closed one \cite{C}. We recall now
the global properties of $M$.

Let $p\in M$ and let $\ga$ be the normal geodesic at $p$.
We know  that $\ga\cap M$ consists of twelve points
$p_1,\dots,p_{12}$ which are vertices of certain dodecagon: 
see Fig.1, where indices are changed from 
[M1, pp. 197--8] and [M2, Lemma 3.2].
\bel \mbox{\rm \cite{M1}} We have the relations
$$
\begin{array}{ll}
D_i(p_1) & =D_{2-i}(p_2)=D_{i+4}(p_3)=D_{4-i}(p_4)
=D_{i+2}(p_5)=D_{6-i}(p_6)\\
D_i(p_j) & =D_i(p_{j+6}),\quad j=1,\dots,6
\end{array}
\label{3.1}
$$
where the equality means ``be parallel to with respect to
the connection of
$S^{7}$", and the indices are modulo 6.
\enl
From these, some relations among $\La_{\al\be}^\ga$'s 
are obtained as follows. 
Denote by $p(t)$ the point on $L_6(p)$ such that $p_1=p(0)$, prametrized by
the center angle where the center means that of a 
circle on a plane. Similarly,  
we denote by $q(t)$ the point on $L_2(p_2)$ parametrized 
from $p_2=q(0)$. Note that $e_6(p_1)$ is parallel with $e_2(p_2)$. 
Extend $e_6$ and $e_2$ as the unit tangent vectors of $p(t)$
and $q(t)$, respectively. Consider the normal geodesic 
$\ga_t$ at  $p(t)$, then 
$q(t)=L_2(p_2)\cap \ga_t$.
Here  $e_{ 3}(p(t))$ is parallel with
$e_{5}(q(t))$. Then we have 
$$
\frac1{\sin\theta_6}\nabla_{\frac{d}{dt}}e_{3}(p(t))
= \frac{\sin\theta_2}{\sin\theta_6}
\frac1{\sin\theta_2}\nabla_{\frac{d}{dt}}e_{5}(q(t)).
$$
Therefore the $D_j$ component of
$(\nabla_{e_{6}}e_{3})(p_{1})$ is  the
$D_{2-j}$ component of
$(\nabla_{e_{2}}e_{5})(p_{2})$ multiplied
by $\sin\theta_2/\sin\theta_6$. We denote such relation by
$$
\La_{{ 6}{3}}^{j}(p_1)
\sim \La_{{2}{5}}^{{2-j}}(p_2),
$$
up to sign. 
A similar argument at every $p_{m}$ implies the global
correspondence among $\La_{\al\be}^\ga$'s:
\bel \mbox{\rm \cite{M1}} For a frame consisting of principal vectors 
around each $p_m$,
we have the correspondence 
$\La_{jk}^{i}(p_m)\sim
\La_{ j'k'}^{ i'}(p_n)$
where $i,j,k$ at $p_m$ correspond to
$i',j',k'$ at $p_n$ in Table 1:
\begin{figure}[htbp]
\begin{center}
\begin{tabular}{c@{\hspace{2em}}c}
\includegraphics[width=6cm]{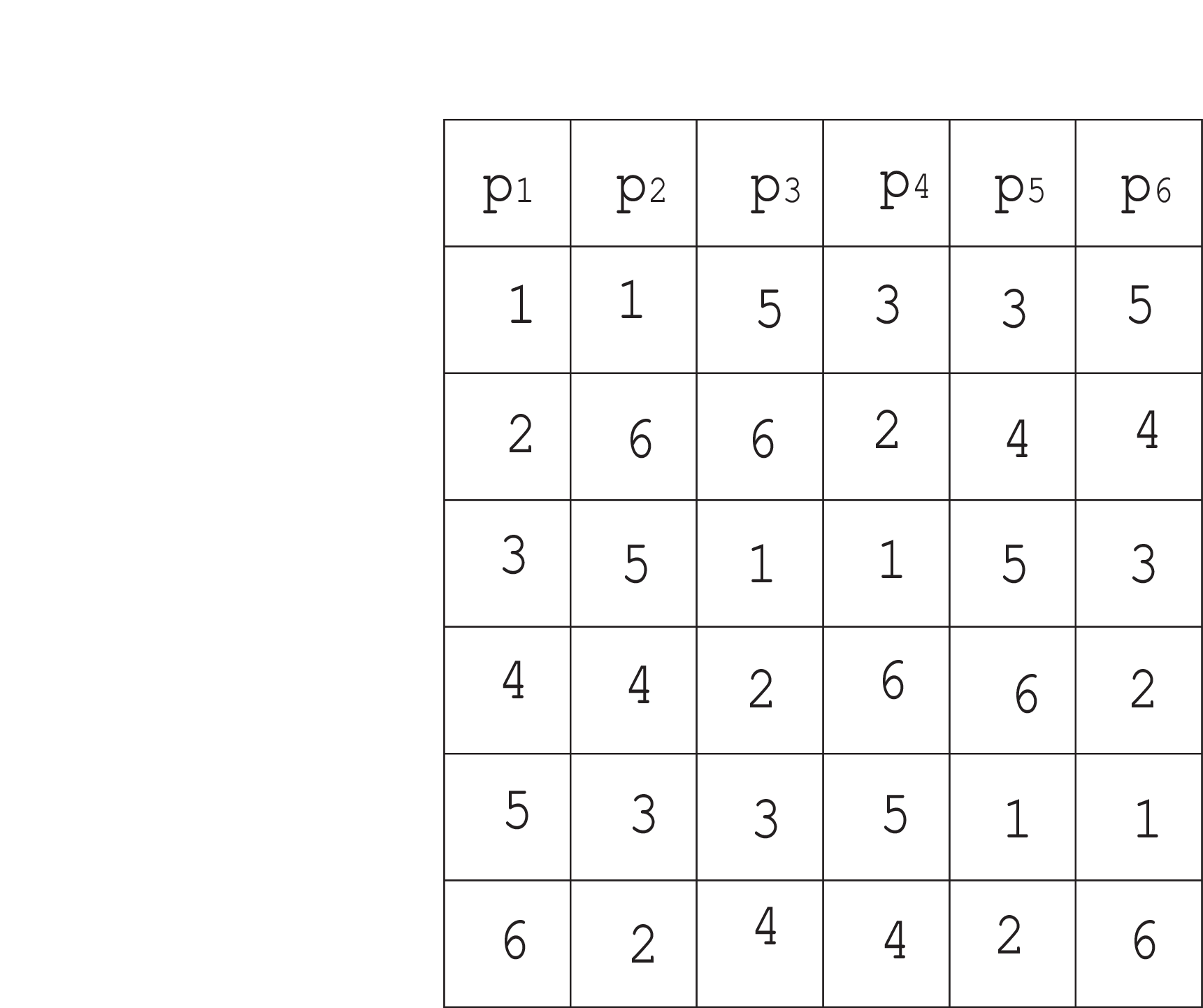}&
\includegraphics[width=5cm]{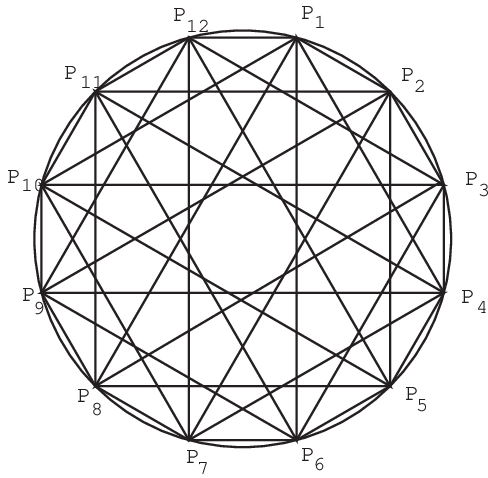}\\
{\small $\qquad\qquad\qquad$Table 1}&{\small Fig.1}
\end{tabular}\end{center}\end{figure}
\label{lmGC}
\enl

%%%%%%%%%%%%%

\section{The kernel of the shape operators}%5
\label{kernel}

For $p\in M$ and $\bar p\in M_+$, let  
$$
E_{\bar p}=\text{span}\{\text{Ker}L(t)\mid t\in [0,2\pi) \}
=\text{span}_{t\in[0,2\pi)}\{e_3(t)\}.
$$
 The following proposition proved in \cite{M1} is crucial.

\bep \mbox{\rm\cite{M1}(Proposition 4.2)} $M$ is 
homogeneous if and only if $\dim E_{\bar p}=1$ for any $p$.
\label{keyp}
\enp
 
Next, recall 
\begin{equation}
\mu_i=\frac{1+\la_i\la_6}{\la_6-\la_i}
=c_1\frac{\la_3-\la_i}{\la_6-\la_i},
\quad c_1=2+\sqrt3.%\label{4.1}
\end{equation}
The second equality follows from
$\la_6=-1/\la_3=-(2+\sqrt3)$.
 Put 
$$
 c_2=\frac1{\sin\theta_6(\la_3-\la_6)}
=-\frac{\sqrt2(\sqrt3+1)}{4},  \quad 
(\sin\theta_6=-\frac{\s2(\s3-1)}{4}).
$$

\bel Take $p\in f_6^{-1}(\bar p)$ and 
identify $T_{\bar p}M_+$ with $\oplus _{j=1}^5D_j(p)$.  
Then we have
\begin{align}
&B_\zeta(e_{ 3})=c_2\nabla_{e_{ 3}}e_6,
\label{43}\\
&B_\eta(\nabla_{e_{ 6}}e_{ 3})=
c_1\nabla_{e_{ 3}}e_{ 6}\label{42}\\
%%0804065%
&B_\zeta(\nabla_{e_{ 6}}e_{ 3})=
c_2\nabla_{e_{ 6}}\nabla_{e_{ 3}}e_{ 6}\label{45}. 
\end{align}
Similar formulas hold for the shape operators 
$C_N$ of $M_-$, if we replace 6 by 1, and 3 by 4. 
\label{lm41}
\enl
\Proof From (\ref{bij}) follows (\ref{43}). 
Using (\ref{2.6}), we have (\ref{42}):
\beq
B_\eta(\nabla_{e_{6}}e_{3})
=\La_{{6}{3}}^{i}
\mu_ie_{i}
=c_1\La_{{6}{3}}^{i}
\frac{\la_3-\la_i}{\la_6-\la_i}e_{i}
=c_1\La_{{3}{6}}^{i}e_{i}
=c_1\nabla_{e_{3}}e_{6}. 
\label{6336}
\enq
Taking the covariant derivative of (\ref{43}) where 
$\nabla_{\frac{d}{dt}}=c_0\nabla_{e_6}$ by (\ref{c0}), 
we obtain 
$$
\bary
c_2\nabla_{e_{ 6}}\nabla_{e_{ 3}}e_{ 6}=\nabla_{e_6}\left(B_\zeta(e_{3})\right)
=-1/c_0B_\eta(e_{3})+B_\zeta(\nabla_{e_{6}}e_{3})
=B_\zeta(\nabla_{e_{6}}e_{3}). 
\eary
$$
\qed
\Remark (\ref{43}) implies that $\dim E_{\bar p}=1$ holds 
 if and only if   $\nabla_{e_{ 6}}e_{ 3}$ vanishes at a 
 point of $f^{-1}(\bar p)$. 
Moreover,  (\ref{42}) implies that 
$\nabla_{e_{ 6}}e_{ 3}$ vanishes if and only if
$\nabla_{e_{ 3}}e_{ 6}$ vanishes. 

\medskip
When $\nabla_{e_6}e_3(p)\not\equiv 0$, 
we have $\dim E_{\bar p}\ge 2$, since 
$e_3(p)$ and $\nabla_{e_6}e_3(p)$ ($\in E_{\bar p}$) 
 are mutually orthogonal.  
We denote $E$ instead of $E_{\bar p}$, 
when it causes no confusion. 
Let $E^\bot$ be the orthogonal complement of $E$ in
$T_{\bar p}M_+$.
Moreover, put  
$$
W=W_{\bar p}=\mbox{\rm span}_{t\in[0,2\pi)}
\{\nabla_{e_3}e_6(t)\}.
$$
where we regard $W$ as a subspace of $T_{\bar p}M_+$ by 
a parallel displacement. 
The following lemmas are significant.

\bel \mbox{\rm\cite{M2} (Lemma 4.2)} $W\subset E^\bot$.
\label{l4.2}
\enl
 
\Proof We can express
$L(t)$ with respect to the basis
$e_i(p)$, $i=1,\dots 5$, as in Lemma \ref{l2.1}, 
$$
L(t)=\begin{pmatrix}
\sqrt3c & s b_{12} & s  b_{13}
&sb_{14} & sb_{15} \\
sb_{12}  & \frac1{\sqrt3}c
 & sb_{23} & sb_{24} & sb_{25} \\
sb_{13} & sb_{23} &  0 & sb_{34}
&  sb_{35} \\
sb_{14} &  sb_{24} & sb_{34}
& -\frac1{\sqrt3}c  &  s  b_{45} \\
sb_{15} &  sb_{25} & sb_{35}
& sb_{45} & -\sqrt3c
\end{pmatrix}, \quad \left\{\begin{array}{ll}
 c=\cos t\\
s=\sin t.
\end{array}\right.
$$
Let 
$e_{ 3}(t)={}^t(u_1(t),
\dots, u_5(t))$
belong to the kernel of $L(t)$. 
Then the third component of
$L(t)(e_{ 3}(t))$ must satisfy
$$
\frac{\sin t}{\sin\theta_6}\frac1{\la_3-\la_6}
\sum_{i=1}^5 
\La_{36}^i(p) u_{i }(t)=0.
$$
Thus we obtain
\begin{equation}
\langle \nabla_{e_3}e_6(p),e_{ 3}(t)\rangle =0
\label{5.1}
\end{equation}
for {\it all} $t$, 
which means 
$\nabla_{e_3}e_6(p)\in E^\bot$.
\qed

By the analyticity and  the definition of $E$ and $W$, 
we can express 
\begin{equation}
\begin{array}{ll}
E &=\text{span}\{e_{ 3}(q),
\nabla_{e_6}^ke_{ 3}(q),\,k=1,2,\dots\}\\
W &=\text{span}\{\nabla_{e_{ 3}}e_6(q),
\nabla_{e_6}^k\nabla_{e_{ 3}}e_6(q),\,k=1,2,\dots\},
\end{array}
\label{anal}
\end{equation}
at any fixed point $q\in L_6$, where $\nabla^k_{e_6}$ means 
$k$-th covariant differential in the direction $e_6$.
Thus we have by Lemma \ref{l4.2}, 
\begin{equation}
%\left\{\begin{array}{ll}
%0=\langle \nabla_{e_3}e_6,
%\nabla_{e_6}^ke_{{3}}\rangle,\\
\langle\nabla_{e_6}^ke_{{3}},\nabla^l_{e_6}\nabla_{e_3}e_{6}\rangle=0,
%\end{array}\right. 
\quad k,l=0,1,2,\cdots.
\label{5.2}
\end{equation}

\bel \mbox{\rm\cite{M2} (Lemma 4.3)} For any $t$,  $L(t)$ maps $E$ onto $W\subset E^\bot $.
\label{14.3}
\enl

\Proof First we show if $L(t)(\nabla_6^ke_{ 3}(p))\in W$ 
holds for any $0\le k\le l$ and $t$, then  
$L(t)(\nabla_6^{l+1}e_{ 3}(p))\in W$ follows. 
In fact, from $L(t)=\cos tB_\eta+\sin tB_\zeta$, 
we have 
$$
L_t(t) =L(t+\pi/2),\quad L_{tt}(t) =-L(t)
$$
Thus in each relation
$$
\bary
L_t(t)(\nabla_6^le_{3}(t))
=c_0\nabla_{e_{ 6}}(L(t)(\nabla_6^le_{ 3}(t)))
-L(t)(c_0\nabla_6^{l+1}e_{ 3}(t)),\\
L_{tt}(t)(\nabla_6^le_{3}(t))
=c_0\nabla_{e_{ 6}}(L_t(t)(\nabla_6^le_{ 3}(t)))
-L_t(t)(c_0\nabla_6^{l+1}e_{ 3}(t)),
\eary
$$
where we use the moving frame $e_3(t)$ and $e_6(t)$, 
the left hand side belongs to $W$ by the assumption,
and so is the first term of the right hand side.
Hence we have 
$L(t)(\nabla_6^{l+1}e_{ 3}(t)),L_t(t)(\nabla_6^{l+1}e_{ 3}(t))\in W$.
Now, we show the lemma by induction. Indeed, 
$L(t)$ maps $D_3(p)$ into  $W$
for all $t$, because $B_{\eta_p}$ and 
$B_{\zeta_p}$ map $D_3(p)$ into $W$ by (\ref{43}),   
and because $L(t)=\cos tB_\eta+\sin tB_\zeta$. 
Moreover,  (\ref{43}) implies 
that this is an onto map.
\qed

\bel \mbox{\rm\cite{M2} (Lemma 4.4)} $\dim E\le 3$. 
\label{l4.4}
\enl
\Proof Take any $p\in f_6^{-1}(\bar p)$. 
Since Ker$B_{\eta_p}=D_3(p)\subset E$, we have 
dim $B_\eta(E)=\dim E-1$.
Because $B_{\eta_p}(E)$ is a subspace of
$E^\bot$, the lemma follows from 
$\R^{5}\cong T_{\bar p}M_+=E\oplus E^\bot$.
\qed

The following is obvious:
\bel As a function of $\bar p\in M_+$,  
dim $E$ is lower-semi-continuous.
\label{lmlsc}
\enl

Let $d=\max_{\bar p\in M_+}\dim E_{\bar p}$. 
We know that $1\le d \le 3$ 
and $M$ is homogeneous when $d=1$. 
At a point $\bar q$ on the focal submanifolds 
$M_-=M_1$, denote 
$F_{\bar q}=\text{span}_{q(t)\in L^1(q)}\{e_4(q(t))\}$.
The argument on $M_+$ holds for $M_-$ if we replace
$E$ by $F$ and pay attention to the change of indices.
Especially, $\dim E=1$ holds on $M_+$ if and only if 
$\dim F=1$ holds on $M_-$, because $\La_{36}^j=0$ holds 
for all $j$ 
if and only if $\La_{14}^j=0$ holds for all $j$,  
by the global correspondense in \S\ref{GP}. 
Note that, however, not everything is symmetric on $M_{\pm}$. 
Indeed, for homogeneous hypersurfaces with six 
principal curvatures, $M_+$ and $M_-$ are {\it not} 
congruent (\S\ref{pre},  \cite{M1}). 

%%%%%%%%%%%%%

 %%%%%%%%%%%%%%%
\section{Description of $E$}
\label{sdescription}

In this section, we discuss what 
happens if we suppose $\dim E\ne 1$. 
Lemma \ref{14.3} suggests that the matrix expression 
of $L(t)$ can be simplified if we use the decomposition 
$T_{\bar p}M_+=E\oplus E^\bot$.

\bel When $\dim E=d$, we can express 
$L=L(t)$ as
$$
L=\begin{pmatrix}
0_d & R \\
{}^tR &S
\end{pmatrix}, 
$$
with respect to the decomposition 
$T_{\bar p}M_+=E\oplus E^\bot$,
where $0_d$ is $d$ by $d$, $R$ is $d$ by $5-d$
and $S$ is $5-d$ by $5-d$ matrices.
The kernel of $L$ is given by 
$$
\begin{pmatrix}X\\0\end{pmatrix}\in E,\quad
{}^tRX=0.
$$
The eigenvectors with respect to $\mu_i(\ne 0)$ 
in (\ref{beta}) are given by 
$$
\begin{pmatrix}\frac1{\mu_i}RY\\Y\end{pmatrix}
$$ 
where $Y\in E^\bot$ is a solution of
\beq
({}^tRR+\mu_iS-\mu_i^2I)Y=0.
\label{eqreduced}
\enq
\label{funda}
\enl

\Proof The first part follows from Lemma \ref{14.3}.
 Let $\begin{pmatrix} X\\Y\end{pmatrix}$ be an 
 eigenvector of $L$ with
respect to $\mu_i$, where $X\in E$ and $Y\in E^\bot$, 
abusing the notation $X=\begin{pmatrix} X\\0\end{pmatrix}$ and 
$Y=\begin{pmatrix} 0\\Y\end{pmatrix}$. Then we have
$$
\begin{pmatrix}
0_d & R \\
{}^tR &S\end{pmatrix}\begin{pmatrix} X\\Y\end{pmatrix}
=\begin{pmatrix} RY \\ {}^tRX+SY\end{pmatrix}
=\mu_i\begin{pmatrix} X\\Y\end{pmatrix},
$$
and hence 
$$
\left\{
\begin{array}{ll}
RY=\mu_iX\\
 {}^tRX+SY=\mu_iY.
\end{array}\right.
$$
For $\mu_3=0$, $Y=0$ and ${}^tRX=0$ hold since the 
kernel belongs to $E$. 
When $\mu_i\ne 0$, multiplying $\mu_i$ to the second
equation and substitute the first
one into it, we obtain (\ref{eqreduced})
Then the eigenvector of $L$ for 
an eigenvalue $\mu_i$ is given by 
$$
\begin{pmatrix}\frac1{\mu_i}RY\\Y\end{pmatrix}.
$$ 
\qed

%%%%%%%%%%%%
\section {Dim $E=2$}
\label{sdimE2}

In this section, we suppose $\dim E=2$ occurs at some point 
$\bar p\in M_+$. Then we have the decomposition 
$T_{\bar p}M_+=E^2\oplus V^2\oplus W^1$
 (the upper indices mean dimensions), where
 $W=B_\eta(E)=B_\zeta(E)$ by Lemma \ref{14.3}.

For a continuous frame $e_3(t)\in D_3(t)$ along $L_6$, 
$D_3(t+\pi)=D_3(t)$ implies $e_3(t+\pi)=\ve e_3(t)$, 
$\ve=\pm1$. Then we have
$\nabla_{e_6}e_3(t+\pi)=\ve \nabla_{e_6}e_3(t)$, and 
it follows 
$$
\bary
\nabla_{e_3}e_6(t+\pi)&=1/c_1L(t+\pi)(\nabla_{e_6}e_3(t+\pi))\\
&=-1/c_1L(t)(\ve \nabla_{e_6}e_3(t))=-\ve\nabla_{e_3}e_6(t).
\eary
$$
Since $\nabla_{e_3}e_6(t)\in W$ never vanishes (Remark 5.3), 
and so has a constant direction, we have $\ve=-1$. 

 In the following, we mean by a continuous frame 
 $e_i(t)$ along $L_6$, a frame on $L_6$ minus a point. 
This is because we may have $e_i(t+2\pi)=-e_i(t)$, 
which occurs as $O(5)$ acts on the shape operator via spin 
action. Fortunately, this does not affect the argument. 

Consider a continuous frame $e_i(t)$ along $L_6$, 
and express $\nabla_{e_6}e_3(t)=\La_{63}^i(t)e_i(t)$. 
Then putting $f(t)=\left(\La_{63}^1(t)\right)^2-\left(\La_{63}^5(t)\right)^2$, we have $f(t+\pi)=-f(t)$ since 
$\nabla_{e_6}e_3(t+\pi)=-\nabla_{e_6}e_3(t)$ and 
$D_i(t+\pi)=D_{6-i}(t)$ holds. 
Thus at some point $p=p(t_0)$ of $L_6$, $f(t_0)=0$ occurs.  
Here by the Gauss equation [3.3], or from 
$$
\bary
0&=\langle \nabla_{e_6}e_3(t),L(t)(\nabla_{e_6}e_3)(t)\rangle \\
&=\s3\{\left(\La_{63}^1(t)\right)^2-\left(\La_{63}^5(t)\right)^2\}+1/\s3\{\left(\La_{63}^2(t)\right)^2-\left(\La_{63}^3(t)\right)^2\},
\eary
$$
we have also 
$\left(\La_{63}^2(t_0)\right)^2-\left(\La_{63}^4(t_0)\right)^2=0$. 
Thus we may put at $p$, 
\beq
\bary
\nabla_{e_6}e_3=x(e_1+e_5)+y(e_2+e_4)\\
\nabla_{e_3}e_6=\s3x(e_1-e_5)+y/\s3(e_2-e_4)
\eary
\label{xy}
\enq
by rechoosing the directions of $e_i=e_i(p)$, 
$i=1,2,4,5$, if necessary. 
Nomalizing the right hand side, we define 
$$
\bary
X_1=\al(e_1+e_5)+\be(e_2+e_4)\in E,\\
Z_1=1/\si\{\s3\al(e_1-e_5)+\be/\s3(e_2-e_4)\}
\in W
\eary
$$
where $\al^2+\be^2=1/2$ and $\si=2(3\al^2+\be^2/3)$, and   
$\nabla_{e_6}e_3=aX_1$ and $\nabla_{e_3}e_6=bZ_1$ hold for 
some $a$ and $b$. Note that $B_\eta(X_1)=\s{\si}Z_1$. 
Since $V$ is orthogonal to $e_3,X_1,Z_1$, 
we have an orthonormal basis of $V$ given by 
$$
\bary
X_2=1/\si\{\be/\s3(e_1-e_5)-\s3\al(e_2-e_4)\},\\
Z_2=\be(e_1+e_5)-\al(e_2+e_4).
\eary
$$
where $B_\eta(X_2)=1/\s{\si}Z_2$ holds. 
Since $V$ is parallel, 
$$
X_2(t)=X_2(0),\quad Z_2(t)=Z_2(0)
$$ 
is an orthonormal frame of $V$ at any $p(t)$. 
Now express $X_2(\pi)=X_2(0)$ and $Z_2(\pi)=Z_2(0)$ via 
basis at $p(\pi)$. Namely, choosing  
$e_i(\pi)=e_i'\in D_i(\pi)=D_{6-i}(0)$ suitably, 
we can express 
$$
\bary
X_2(\pi)&=1/\si\{\be'/\s3(e'_1-e'_5)-\s3\al'(e'_2-e'_4)\}\\
&=1/\si\{\be/\s3(e_1-e_5)-\s3\al(e_2-e_4)\},\\
Z_2(\pi)&=\be'(e'_1+e'_5)-\al'(e'_2+e'_4)\\
&=\be(e_1+e_5)-\al(e_2+e_4),
\eary
$$
because $D_1(\pi)\oplus D_5(\pi)=D_1(0)\oplus D_5(0)$, 
and  $D_2(\pi)\oplus D_4(\pi)=D_2(0)\oplus D_4(0)$ hold, 
and hence
$|\al'|=|\al|$, $|\be'|=|\be|$, and $\si=\si(\pi)=\si(0)$ 
follow. Thus we obtain 
$$
\left\{
\bary
\be'(e'_1-e'_5)= \be(e_1-e_5),\quad\\
\be'(e'_1+e'_5)= \be(e_1+e_5), 
\eary
\right.\quad 
\left\{\bary
\al'(e'_2-e'_4)=  \al(e_2-e_4), \\ 
\al'(e'_2+e'_4)= \al(e_2+e_4),  
\eary
\right.
$$
and from $D_i(\pi)=D_{6-i}(0)$, it follows 
$$
\left\{\bary
\be' e'_1=- \be e_5,\\
-\be' e'_5= \be e_1,\\
\be' e'_1= \be e_5,\\
\be' e'_5= \be e_1,
\eary\right.\quad
\left\{\bary
\al' e'_2=- \al e_4,\\
-\al' e'_4= \al e_2,\\
\al' e'_2= \al e_4,\\
\al' e'_4= \al e_2.
\eary\right.
$$
However then, we have $\al=\be=0$, a contradiction. 

Thus we conclude:
\bep
$\dim E=2$ does not occur at any point of $M_+$.
\label{pdimE2}
\enp 

%%%%%%%%%%%
\section {Dim $E=3$}
\label{sdimE3}

By the previous proposition, $\dim E=3$ occurs on $M_+$ if 
$\dim E>1$.  

\bep When $\dim E=3$, at any point $p$ of $L_6$, $E$ and 
$E^\bot$ are expressed via $e_i=e_i(p)$ as
$$
\bary
E=\text{span}\{e_3,\al(e_1+e_5)+\be(e_2+e_4), 
\frac{\be}{\s3}(e_1-e_5)-\s3\al(e_2-e_4)\}\\
 E^\bot=\text{span}\{\s3\al(e_1-e_5)+\frac{\be}{\s3}(e_2-e_4),
\be(e_1+e_5)-\al(e_2+e_4)\},
\eary
$$
for suitable $\al,\be$ satisfying $\al^2+\be^2\ne0$. 
\label{peachcase}
\enp
%%%%%080410%%%%%
\Proof 
Since $e_3,e_1+e_5,e_2+e_4,e_1-e_5,e_2-e_4$ generate a frame 
of $TM_+$, we can choose $X_1,X_2\in E$ as 
$$
\bary
X_1=\al(e_1+e_5)+\be(e_2+e_4)+\ga(e_1-e_5)\\
X_2=x(e_1+e_5)+y(e_2+e_4)+z(e_1-e_5)+w(e_2-e_4).
\eary
$$
Then $Z_i=B_\eta(X_i)\in E^\bot$ are given by 
$$
\bary
Z_1=\s3\al(e_1-e_5)+\frac1{\s3}\be(e_2-e_4)+\s3\ga(e_1+e_5)\\
Z_2=\s3x(e_1-e_5)+\frac1{\s3}y(e_2-e_4)+\s3z(e_1+e_5)+\frac1{\s3}w(e_2+e_4).
\eary
$$
Because $0=\langle X_1,Z_1\rangle=2\s3\al \ga$, 
changing the sign of $e_5$, if necessary, we may assume 
$\ga=0$, i.e., 
\beq
\bary
X_1=\al(e_1+e_5)+\be(e_2+e_4)\in E\\
Z_1=\s3\al(e_1-e_5)+\frac{\be}{\s3}(e_2-e_4)\in E^\bot. 
\eary
\label{XX12}
\enq
Next from 
$0=\langle X_1,Z_2\rangle=\s3\al z+\frac{\be w}{\s3}$,  
and $0=\langle X_2,Z_2\rangle=2(\s3x z+\frac1{\s3}y w)$, 
$\al y-\be x=0$ holds unless $z=w=0$, and then 
$x(e_1+e_5)+y(e_2+e_4)$ is proportional to $X_1$. Thus 
we may rechoose 
\beq
\bary
X_2=z(e_1-e_5)+w(e_2-e_4)
=\frac{\be}{\s3}(e_1-e_5)-\s3\al(e_2-e_4)\in E,
\eary
\label{X21}
\enq
and  
\beq
Z_2=\be(e_1+e_5)-\al(e_2+e_4)\in E^\bot.
\label{Z2}
\enq
When $z=w=0$, we have span$\{X_1,X_2\}=\text{span}\{e_1+e_5,e_2+e_4\}$  
and span$\{Z_1,Z_2\}=\text{span}\{e_1-e_5,e_2-e_4\}$.  
Here, in order to fit in the expression (\ref{X21}) and 
(\ref{Z2}), we change the sign of $e_4$, and may consider
\beq
X_2=e_2-e_4,\quad Z_2=e_2+e_4,
\enq
corresponding to $\be=0$. 
\qed
%%%%%%%%%
Note that $X_1,X_2,Z_1,Z_2$ are mutually orthogonal. 
Then the orthonormal frames of $E$ and $E^\bot$ are 
given respectively, by
\beq
\bary
e_3,\quad X_1=\al(e_1+e_5)+\be(e_2+e_4)\\
X_2=\dfrac{1}{\s{\si}}\left(\frac{\be}{\s3}(e_1-e_5)-\s3\al(e_2-e_4)\right)
\eary
\label{basisE}
\enq
and
\beq
\bary
Z_1=\dfrac{1}{\s{\si}}\left(\s3\al(e_1-e_5)+\frac{\be}{\s3}(e_2-e_4)\right)\\
Z_2=\be(e_1+e_5)-\al(e_2+e_4),
\eary
\label{basisEbot}
\enq
where we put 
\beq
\al^2+\be^2=1/2,\quad \si=2(3\al^2+\be^2/3).
\label{albe}
\enq 
Consider an arc $c$ of $L_6$ containing $p=p(0)$ and $p(\pi)$. 
Since $X_1,X_2$ are given at each point of $L_6$ by (\ref{XX12}), (\ref{X21}) and  
(\ref{Z2}), using a continuous frame $e_i(t)$ and 
a continuous function $\al(t), \be(t)$ along $c$, 
we have a continuous frame $e_3(t)$, 
$X_1(t)$ and $X_2(t)$ of $E$, and $Z_1(t)$ and $Z_2(t)$ of 
$E^\bot$ along $c$. 
With respect to this moving frame, we can express
\beq
L(t)=B_{\eta_t}=\bpm 0&0&0&0&0\\
0&0&0&\s{\si(t)}&0\\
0&0&0&0&1/\s{\si(t)}\\
0&\s{\si(t)}&0&0&u(t)\\
0&0&1/\s{\si(t)}&u(t)&0\epm 
\label{Beta}
\enq
for $\eta_t=\eta_{p(t)}$DIn fact, from 
$L(t)(e_i(t))=\mu_ie_i(t)$, we know 
$L(t)(X_1(t))=\s{\sigma}Z_1(t)$ and 
$L(t)(X_2(t))=1/\s{\sigma}Z_2(t)$. 
Moreover, it is easy to see 
$\langle L(t)(Z_i(t)), Z_i(t)\rangle=0$. 
Then putting $u(t)=\langle  L(t)(Z_1(t)), Z_2(t)\rangle$, 
we have (\ref{Beta}). 
Note that $\si(t)+1/{\si(t)}+u(t)^2={10}/{3}$ follows from  
$\|L(t)\|=\frac{20}{3}$. Moreover, by using the notation 
in \S\ref{sdescription}, (\ref{Beta}) implies that 
$T(t)={}^tR(t)R(t)$ 
has eigenvalues ${\si(t)},1/{\si(t)}$ with eigenvectors 
$Z_1(t),Z_2(t)\in E^\bot$, respectively. 
Note that even if ${\si(t)}=1/{\si(t)}$ holds, 
$Z_1(t)$ and $Z_2(t)$ (thus, $X_1(t)$ and $X_2(t)$) are 
continuously chosen so that the $S(t)$ part in (\ref{Beta}) be 
described as above where $u(t)^2=4/3\ne 0$.  

Next, we show: 
%%%%%%%0741%%%%%
\bep $\si(t)$ is constant and takes the values 1, 1/3 or 3. 
\label{pnot}
\enp
\Proof We have $L(\pi)=-L(0)$ from 
$L(t)=\cos t B_\eta+\sin tB_\zeta$, and 
$T(\pi)=T(0)$ from $T(t)={}^tR(t)R(t)$. This implies  
$\si=\si(\pi)=\si(0)$.   
When $\si(t)$ is not identically 1,  
we may consider $\si\ne 1$, 
and as an eigenvector of $T(0)$ for $\si$,
$Z_1(\pi)$ is parallel to $Z_1(0)$. Then from 
$$
\left\{\bary
L(\pi)(X_1(\pi))=\s{\si}Z_1(\pi),\\
L(0)(X_1(0))=\s{\si}Z_1(0),
\eary\right.
$$
 we have
$$
 X_1(\pi)=\ve X_1(0), \quad Z_1(\pi)=-\ve Z_1(0),\quad 
 \ve=\pm1.   
 $$ 
Similarly  from 
 $$
\left\{\bary
L(\pi)(X_2(\pi))=1/\s{\si}Z_2(\pi),\\
L(0)(X_2(0))=1/\s{\si}Z_2(0),
\eary\right.
$$
 we have, unless $\al\be\not\equiv 0$, 
$$
 X_2(\pi)=-\ve X_2(0),\quad Z_2(\pi)=\ve Z_2(0), 
\label{X2}
 $$
where we use $e_i(\pi)\in D_{6-i}(0)$ by the global 
correspondence in (\ref{basisE}) and (\ref{basisEbot}). 
However, since $E^\bot$ is parallel along $L_6$,  
and the pair $Z_1(t),Z_2(t)$ is a continuous orthonormal 
frame of $E^\bot$ by the remark before the proposition, 
this contradicts the fact that a continuous frame 
preserves the orientation. 
Therefore, only the cases $\si\equiv 1,1/3,3$ remain. 
\qed  

\section {Final result}
\label{sfinal}

\bep When $\dim E=3$, $\si\equiv 1$ does not occur.
\label{psi1}
\enp
\Proof In this case, $3\al^2=\be^2$ follows from (\ref{albe}), 
and hence by a suitable choice of directions of $e_i$'s, we have  
$$
\bary
E=\text{span}\{e_3, e_1+\s3e_4,\s3e_2+e_5\}\\
E^\bot=\text{span}\{\s3e_1-e_4,e_2-\s3e_5\}.
\eary
$$
Since $B_\zeta$ maps $E$ onto $E^\bot$, 
$b_{14}=b_{25}=0$ follows, i.e., 
$\La_{16}^4=\La_{26}^5=0$ holds. These imply 
$\La_{63}^2=\La_{63}^4=0$ by the global correspondence. 
However, since $\nabla_{e_6}e_3$ is a combination of 
$e_1+\s3e_4$ and $\s3e_2+e_5$, this implies 
$\nabla_{e_6}e_3=0$,  a contradiction. 
\qed

In the last possible case, we have by Proposition 
\ref{peachcase}, 
$$
E=\text{span}\{e_3, e_1+e_5,e_2-e_4\},\quad
E^\bot=\text{span}\{e_1-e_5,e_2+e_4\},
$$
and this holds everywhere by a continuous choice of 
$e_i$'s. 
Since $E$ is mapped onto $E^\bot$ by $B_\zeta=(b_{ij})$, 
we have  
\beq
b_{15}=b_{24}=0,\quad b_{12}+b_{25}=b_{14}+b_{45}.
\label{b1524}
\enq
On the other hand, for another focal submanifold $M_-$, 
the remaining possible case is also this case when dim $F=3$. 
(For the definition of $F$, see the end of \S\ref{kernel}.) Because 
$\nabla_{e_3}e_6(p)\sim \nabla_{e_1}e_4(q)\in E^\bot\cap F$, 
where $p=p_1$ and $q=p_3$ in Fig 1, 
identifying the vectors at $q$ with those at $p$ 
as in Table 1, we may consider
$$
\bary
F&=\{e_4(q), e_5(q)-e_3(q), e_6(q)+e_2(q)\}\\
&=\{e_6(p), e_1(p)-e_5(p),e_2(p)+e_4(p)\},\\
F^\bot&=\{e_5(q)+e_3(q), e_6(q)-e_2(q)\}\\
&=\{e_1(p)+e_5(p),e_2(p)-e_4(p)\}.
\eary
$$
Here, some signature might be opposite, which does not matter.
The importance is 
$$
c_{35}=c_{26}=0
$$
holds where $c_{ij}=
\frac{1}{\sin\theta_1(\la_i-\la_1)}\La_{i1}^j$ 
is the components of the shape operator $C_\zeta$ of $M_-$ 
for $\zeta=e_1$ (see Lemma \ref{l2.1}). 
Then the latter implies $b_{12}=0$, and by the global correspondence, we have $b_{45}=0$, 
and hence it follows from (\ref{b1524}), 
$$
b_{14}=b_{25}.
$$
Next from the Gauss equation [1.2] in \S\ref{sisospectral}, $b_{13}b_{32}=0$ 
follows. When $b_{13}=0$, [1.1] implies $b_{14}^2=2$, 
and hence $b_{25}^2=2$, but this contradicts [2.2]. 
Thus we have $b_{23}=0$. 
Since this holds identically by the analyticity, 
$b_{14}=b_{25}=0$ follows from the global 
correspondence, 
and the second row of $B_\zeta$ vanishes, 
contradicts [2.2]. Therefore we obtain:
\bep
$\dim E=3$ does not occur.
\label{piii}
\enp

\medskip
Finally, the kernel of the shape operators of the focal submanifolds of isoparametric 
hypersurfaces with $(g,m)=(6,1)$ is independent of the 
normal directions, 
and by Proposition 4.2 of \cite{M1}, we obtain:
\bet \mbox{\rm \cite{DN}} Isoparametric hypersurfaces with $(g,m)=(6,1)$ 
are homogeneous.
\label{tdorf}
\ent

%%%%%%%%%%%%

\begin{flushleft}

Mathematical Institute,\\
Graduate School of Sciences\\
Tohoku University\\
Aoba-ku, Sendai, 980-8578/JAPAN\\
\textit{E-mail Address}: r-miyaok@math.tohoku.ac.jp

\end{flushleft}
\end{document}